\documentclass{amsart}
\usepackage[cp1251]{inputenc}
\usepackage[english]{babel}
\usepackage{amsmath}
\usepackage{amssymb}
\usepackage{amsfonts}

\begin{document}
\thispagestyle{empty}

 \title[ON A NEW FORMULAS FOR A DIRECT AND INVERSE CAUCHY PROBLEMS ]{  ON A NEW FORMULAS FOR A DIRECT AND INVERSE CAUCHY PROBLEMS OF HEAT EQUATION
 }%

\author{N.Yaremko, O.Yaremko}%Указываем авторов

\address{Natalia Yaremko, Oleg Yaremko
\newline\hphantom{iii}Penza State University,% Место работы
\newline\hphantom{iii}str. Lermontov, 37, % Адрес (улица, дом, строение и т.п.)
\newline\hphantom{iii} 440038, Penza, Russia}%  Адрес (почтовый индекс, город, страна)
\email{yaremki@mail.ru}% Ваш электронный адрес для переписки

\maketitle {\small

\begin{quote}
\noindent{\bf Abstract. } In this paper a solution of the direct Cauchy problems
for heat equation is founded in the Hermite polynomial series form. A
well-known classical solution of direct problem is represented in the
Poisson integral form. The author shows the formulas for the solution of the
inverse Cauchy problems have a symmetry with respect to the formulas for the
corresponding direct problems. The obtained solution formulas for the
inverse problems can serve as a basis for regularizing computational
algorithms while well-known classical formula for the solution of inverse
problem did not possess such properties and can't be a basis for
regularizing computational algorithms.
\medskip

\noindent{\bf Keywords:} heat equation, direct/inverse Cauchy problem, well-posed
/ill-posed problem, Hermite polynomials, Poisson integral.

\emph{{Mathematics Subject Classification 2010}:\\{65Rxx		Integral equations, integral transforms; 12E10	Special polynomials}.}
\end{quote} }

\section{Introduction.}
 In this paper a direct and inverse Cauchy problems
for the heat equation are solved at Cartesian and polar coordinates.
The inverse Cauchy problem for the heat equation consists of to
reconstruct a priori unknown initial condition of the dynamic system from
its known final condition. In 1939 French mathematician Jacques Hadamard
defined the problem is called well-posed if a solution exists, the solution
is unique, the solution's behavior hardly changes when there's a slight
change in the initial condition. The problems are called ill-posed or not
well-posed if at least one of these three conditions is not fulfilled. The
most often, the third condition so called the stability condition of
solution is violated for ill-posed problems. In this case there is a
paradoxical situation: the problem is mathematically defined but it's
solution cannot be obtained by conventional methods. In mathematics the vast
majority of inverse problems is not well-posed: small perturbations of the
initial data (observations) may correspond to an arbitrarily large
perturbations of the solution. A classic example of ill-posed problem is the
inverse Cauchy problem (retrospective problem) for heat equation. The direct
Cauchy problem for the heat equation most frequently is well-posed.

In this work the solution for direct Cauchy problems found in the form of
the Hermite polynomial series. A well-known classical solution for the
direct Cauchy problem is represented in the Poisson integral form. We
noticed the solution of inverse Cauchy problem possess the symmetry with
respect to the solution of the direct Cauchy problem. The formulas obtained
in this paper for the inverse problem can serve as a basis for regularizing
 computational algorithm. Previously known classical formula for
the solution of the inverse problem did not have such properties that's why
classical formula can't serve as a basis for convergent algorithm.

The main result for the direct problems is the formulas (\ref{eq50}), (\ref{eq51}), (\ref{eqf}), (5),(6), (8).The main result for the inverse problems is the formulas (\ref{eq4}),
(\ref{eq5}), (\ref{eq6}), (26), (27), (28).

\section{The Cauchy problem (direct problem) for the heat equation.}

A solution $u\left( {\tau ,x} \right)$ of the Cauchy problem for an infinite
bar with the initial thermal field $f\left( x \right)$ we will get in the
Hermite polynomial series form. For this we use the well-known analytic
solution $u\left( {\tau ,x} \right)$ in the Fourier integral form, [1]:
\[
u\left( {\tau ,x} \right)=\frac{1}{2\pi }\int_{-\infty }^\infty {e^{-\lambda
^2t}} e^{i\lambda x}\left( {\int_{-\infty }^\infty {e^{-i\lambda \xi
}\,f\left( \xi \right)d\xi \:} } \right)d\lambda ,
\]
where $f\left( x \right)$ - the initial thermal field, \newline
$u\left( {\tau ,x} \right)$ - thermal field at the time $\tau $ and at the
point $x$.

Write the last equality in the form
\begin{equation}
\label{eqQ}
u\left( {\tau ,x} \right)=\frac{1}{2\pi }\int_{-\infty }^\infty {e^{-\lambda
^2(\tau +\beta )}} e^{i\lambda x}\left( {\int_{-\infty }^\infty {e^{\lambda
^2\beta }e^{-i\lambda \xi }\,f\left( \xi \right)d\xi \:} } \right)d\lambda
\end{equation}

where $\beta >0$.

The function $e^{\lambda ^2\tau +i\lambda x}$ is a generating function for
the Hermite polynomials, [1], this means that
\begin{equation}
\label{eq4}
~e^{\lambda ^2\beta -i\lambda \xi }=\sum\limits_{j=0}^\infty
{\frac{(-i\lambda )^j\beta ^{\frac{j}{2}}}{j!}} H_j \left( {\frac{\xi
}{2\sqrt \beta }} \right) ,
\end{equation}
where
\begin{equation}
\label{eq1}
H_j (z)=(-1)^je^{z^2}\frac{d^j}{dz^j}\left( {e^{-z^2}} \right)
\end{equation}
- the Hermite polynomials.

In accordance with (2) the formula (\ref{eqQ}) takes the form
\[
u\left( {\tau ,x} \right)=\frac{1}{2\pi }\int_{-\infty }^\infty
{e^{-\lambda ^2(\tau +\beta )}} e^{i\lambda x}\sum\limits_{j=0}^\infty
{\frac{(-i\lambda )^j}{j!}} \beta ^{\frac{j}{2}}\int_{-\infty }^\infty {H_j
} \left( {\frac{\xi }{2\sqrt \beta }} \right)f\left( \xi \right)d\xi
d\lambda
\]

If we may change the order of integration, then compute the inner integral
on the variable $\lambda .$ We use for this Poisson integral, [13]. We get
\begin{equation}
\label{eq5}
\frac{1}{2\pi }\int_{-\infty }^\infty {e^{-\lambda ^2(\tau +\beta )}}
e^{i\lambda x}d\lambda =\frac{e^{-\frac{x^2}{4(\tau +\beta )}}}{2\sqrt {\pi
(\tau +\beta )} }. ~
\end{equation}
Differentiating both sides of ( 4 ) with respect to $x$ we obtain the value
of the desired integral
\[
\frac{1}{2\pi }\int_{-\infty }^\infty ( -i\lambda )^je^{-\lambda ^2(\tau
+\beta )}e^{i\lambda x}d\lambda =(-1)^j\frac{d^j}{dx^j}\left[
{\frac{e^{-\frac{x^2}{4(\tau +\beta )}}}{2\sqrt {\pi (\tau +\beta )} }}
\right].
\]
On the basis of formula (3) for the Hermite polynomials the last formula can
be written
\[
\frac{1}{2\pi }\int_{-\infty }^\infty ( -i\lambda )^je^{-\lambda ^2(\tau
+\beta )}e^{i\lambda x}d\lambda =\frac{e^{- \frac{x^2}{4(\tau +\beta
)}}}{2\sqrt {\pi (\tau +\beta )} }\frac{1}{(2\sqrt {\tau +\beta } )^j}H_j
\left( {\frac{x}{2\sqrt {\tau +\beta } }} \right).
\]

Finally, we obtain an analytic representation of the thermal field at the
time $\tau $ and at the point $x$
\begin{equation}
\label{eq50}
u(\tau ,x)=\frac{e^{-\frac{x^2}{4(\tau +\beta )}}}{2\sqrt
{\pi (\tau +\beta )} }\sum\limits_{j=0}^\infty {\frac{1}{(2\sqrt {\tau
+\beta } )^j}} H_j \left( {\frac{x}{2\sqrt {\tau +\beta } }}
\right)\frac{\beta ^{\frac{j}{2}}}{j!}f_j ,
\end{equation}
where
\[
f_j =\int_{-\infty }^\infty {H_j } \left( {\frac{\xi }{2\sqrt \beta }}
\right)f\left( \xi \right)d\xi .
\]
\textbf{Remark}. By interchanging the variables $x$ and $\xi $ at the formula
(\ref{eqQ}) we get another version of the thermal field:

\begin{equation}
\label{eq51}
u(\tau ,x)=\frac{e^{-\frac{x^2}{4\tau }}}{2\sqrt {\pi \tau }
}\sum\limits_{j=0}^\infty {\frac{1}{(2\sqrt \tau )^j}} H_j \left(
{\frac{x}{2\sqrt \tau }} \right)\frac{\left( {\tau +\beta }
\right)^{\frac{j}{2}}}{j!}f_j ,
\end{equation}

Where
\[
f_j =\int_{-\infty }^\infty {H_j } \left( {\frac{\xi }{2\sqrt {\tau +\beta }
}} \right)f\left( \xi \right)d\xi .
\]
Now we get the third new formula. To do this, the formula (\ref{eqQ}) can be written
in the form
\[
u\left( {\tau ,x} \right)=\frac{1}{2\pi }\int_{-\infty }^\infty
{e^{-\lambda ^2(\tau +\beta )}} \left( {\int_{-\infty }^\infty {e^{\lambda
^2\beta }e^{-i\lambda \left( {x-\xi } \right)}\,f\left( \xi \right)d\xi \:}
} \right)d\lambda ,
\]

here $\beta >0$.

In view of (2) we get
\[
u\left( {\tau ,x} \right)=\frac{1}{2\pi }\int_{-\infty }^\infty
{e^{-\lambda ^2(\tau +\beta )}} \sum\limits_{j=0}^\infty {\frac{(-i\lambda
)^j}{j!}} \beta ^{\frac{j}{2}}\int_{-\infty }^\infty {H_j } \left(
{\frac{x-\xi }{2\sqrt \beta }} \right)f\left( \xi \right)d\xi d\lambda  .
\]

To simplify last formula change the order of integration and compute the
inner integral on the variable $\lambda $,

substitute $x=0$ in (4). We obtain
\begin{equation}
\label{eqb}
\frac{1}{2\pi }\int_{-\infty }^\infty ( -i\lambda )^je^{-\lambda
^2(\tau +\beta )}d\lambda =\frac{1}{(2\sqrt {\tau +\beta } )^{j+1}}H_j
\left( 0 \right).
\end{equation}
Taking into account the well-known formula from [1]
\[
H_{2j} \left( 0 \right)=\frac{\left( {-1} \right)^j(2j)!}{2^jj!},H_{2j+1}
\left( 0 \right)=0;j=0,1,2,...
\]

finally we get an analytical representation of the thermal field
\begin{equation}
\label{eqf}
u(\tau ,x)=\sum\limits_{j=0}^\infty {\frac{1}{(2\sqrt {\tau +\beta }
)^{2j+1}}} \frac{\left( {-1} \right)^j\beta ^j}{2^jj!}f_{2j} ,
\end{equation}

Where
\[
f_{2j} =\int_{-\infty }^\infty {H_{2j} } \left( {\frac{x-\xi }{2\sqrt \beta
}} \right)f\left( \xi \right)d\xi .
\]
\subsection{The inverse Cauchy problem for the heat
equation}.

The inverse heat equation problem for an infinite bar is to find the unknown
initial distribution $f\left( x \right)$ of thermal field by the known
temperature field $u\left( {\tau ,x} \right)$ [7],[8],[13]. This inverse heat equation
problem leads to the solving of first - type Fredholm integral equation:
\begin{equation}
\label{eq2}
\int_{-\infty }^\infty {\frac{1}{2\sqrt {\pi \tau } }\exp \left(
{-\frac{\left( {x-\xi } \right)^2}{4\tau }} \right)f\left( \xi \right)d\xi }
=u\left( {\tau ,x} \right),
\end{equation}
The left-hand side of equation (\ref{eq2}) is the Poisson integral, [15]. As it is
shown in [1],[2] the solution of equation (\ref{eq2}) may be given on the form:
\begin{equation}
\label{eqa}
f\left( x \right)=\frac{1}{\sqrt \pi }\sum\nolimits_{j=0}^\infty
{\frac{u^{\left( j \right)}\left( 0 \right)\:}{\left( {2\sqrt \tau }
\right)^{n+\mbox{1}}j\,!}} H_j \left( {\frac{x}{2\sqrt \tau }} \right),
\end{equation}

where $H_j ( z)-$ Hermite polynomials, (3).

Formula (\ref{eqa}) contains a derivatives of arbitrarily high order so the formula
(\ref{eqa}) can't serve as a basis for the regularizing  computational
algorithm. Consequently, it is actual to find the new formulas without
derivatives for the solving of the equation (\ref{eq2}).

~ As in chapter-1, we obtain three new formulas.

We get the~ solution of equation (\ref{eq2}) by the Fourier transform integral
method from [1],[9],[11],[12]

\[
f\left( x \right)=\frac{1}{2\pi }\int_{-\infty }^\infty {e^{\lambda ^2\tau
}} e^{i\lambda x}\left( {\int_{-\infty }^\infty {e^{-i\lambda \xi }\,u\left(
{\tau ,\xi } \right)d\xi \:} } \right)d\lambda
\]
If $\beta >0$, then the last formula takes the form
\begin{equation}
\label{eq3}
f\left( x \right)=\frac{1}{2\pi }\int_{-\infty }^\infty {e^{-\lambda ^2\beta
}} e^{\lambda ^2(\tau +\beta )}e^{i\lambda x}\left( {\int_{-\infty }^\infty
{e^{-i\lambda \xi }\,u\left( {\tau ,\xi } \right)d\xi \:} } \right)d\lambda
,
\end{equation}
Because of the formula (2) we get

\[
f\left( x \right)=\frac{1}{2\pi }\int_{-\infty }^\infty {e^{-\lambda ^2\beta
}} e^{i\lambda x}\sum\limits_{j=0}^\infty {\frac{(-i\lambda )^j}{j!}} (\tau
+\beta )^{\frac{j}{2}}\int_{-\infty }^\infty {H_j } \left( {\frac{\xi
}{2\sqrt {\tau +\beta } }} \right)u\left( {\tau ,\xi } \right)d\xi d\lambda
\]
We will change the order of integration and calculate the inner integral on
the variable$\lambda $. The Poisson integral from [15] is used in this
calculations. We get
\[
\frac{1}{2\pi }\int_{-\infty }^\infty {e^{-\lambda ^2\beta }} e^{i\lambda
x}d\lambda =\frac{e^{-\frac{x^2}{4\beta }}}{2\sqrt {\pi \beta } }.
\]
We calculate the value of an integral by $j$ time differentiating on the
variable $x$:
\[
\frac{1}{2\pi }\int_{-\infty }^\infty ( -i\lambda )^je^{-\lambda ^2\beta
}e^{i\lambda x}d\lambda =(-1)^j\frac{d^j}{dx^j}\left[
{\frac{e^{-\frac{x^2}{4\beta }}}{2\sqrt {\pi \beta } }} \right].
\]
On the basis of the formula (3) we can write
\[
\frac{1}{2\pi }\int_{-\infty }^\infty ( -i\lambda )^je^{-\lambda ^2\beta
}e^{i\lambda x}d\lambda =\frac{e^{-\frac{x^2}{4\beta }}}{2\sqrt {\pi \beta }
}\frac{1}{(2\sqrt \beta )^j}H_j \left( {\frac{x}{2\sqrt {\beta } }}
\right).
\]
Finally, first new formula for the initial distribution of the thermal field
takes the form
\begin{equation}
\label{eq4}
f(x)=\frac{e^{-\frac{x^2}{4\beta }}}{2\sqrt {\pi \beta }
}\sum\limits_{j=0}^\infty {\frac{1}{(2\sqrt \beta )^j}} H_j \left(
{\frac{x}{2\sqrt \beta }} \right)\frac{(\tau +\beta )^{\frac{j}{2}}}{j!}u_j
,
\end{equation}
where
\[
u_j =\int_{-\infty }^\infty {H_j } \left( {\frac{\xi }{2\sqrt {\tau +\beta }
}} \right)u\left( {\tau ,\xi } \right)d\xi .
\]
\textbf{Remark}. By interchanging the variables $x$ and $\xi $ at the formula
(\ref{eq4}) we get a second new formula for the solution of the equation (\ref{eq2}):
\begin{equation}
\label{eq5}
f(x)=\sum\limits_{j=0}^\infty {\frac{1}{(2\sqrt \beta )^j}} H_j \left(
{\frac{x}{2\sqrt {\tau +\beta } }} \right)\frac{(\tau +\beta
)^{\frac{j}{2}}}{j!}u_j ,
\end{equation}
where
\[
u_j =\int_{-\infty }^\infty {\frac{e^{-\frac{\xi ^2}{4\beta }}}{2\sqrt {\pi
\beta } }} H_j \left( {\frac{\xi }{2\sqrt \beta }} \right)u\left( {\tau ,\xi
} \right)d\xi .
\]
Finally we prove the third new formula for the inverse Cauchy problem
solving.

We use the integral representation of the solution (\ref{eq3}) which

can be written as
\[
f\left( x \right)=\frac{1}{2\pi }\int_{-\infty }^\infty {e^{-\lambda ^2\beta
}} e^{\lambda ^2(\tau +\beta )}\left( {\int_{-\infty }^\infty {e^{i\lambda
\left( {x-\xi } \right)}\,u\left( {\tau ,\xi } \right)d\xi \:} }
\right)d\lambda ,
\]
where $\beta >0$.

Because of the formula (\ref{eq4}) we get
\[
f\left( x \right)=\frac{1}{2\pi }\int_{-\infty }^\infty {e^{-\lambda ^2\beta
}} \sum\limits_{j=0}^\infty {\frac{(-i\lambda )^j}{j!}} (\tau +\beta
)^{\frac{j}{2}}\int_{-\infty }^\infty {H_j } \left( {\frac{x-\xi }{2\sqrt
{\tau +\beta } }} \right)u\left( {\tau ,\xi } \right)d\xi d\lambda
\]
If we use the formula (\ref{eqb}), then the initial distribution of the thermal
field takes the form
\[
f(x)=\frac{1}{2\sqrt {\pi \beta } }\sum\limits_{j=0}^\infty
{\frac{1}{(2\sqrt \beta )^j}} H_j \left( 0 \right)\frac{(\tau +\beta
)^{\frac{j}{2}}}{j!}u_j ,
\]
where
\[
u_j =\int_{-\infty }^\infty {H_j } \left( {\frac{x-\xi }{2\sqrt {\tau +\beta
} }} \right)u\left( {\tau ,\xi } \right)d\xi .
\]
Because of the formula for $H_j \left( 0 \right)$ as a result we get
\begin{equation}
\label{eq6}
f(x)=\frac{1}{\sqrt \pi }\sum\limits_{j=0}^\infty {\frac{1}{(2\sqrt \beta
)^{2j+1}}\frac{\left( {-1} \right)^j(\tau +\beta )^j}{2^j(2j)!}} u_{2j} ,
\end{equation}
Where
\[
u_j =\int_{-\infty }^\infty {H_j } \left( {\frac{x-\xi }{2\sqrt {\tau +\beta
} }} \right)u\left( {\tau ,\xi } \right)d\xi .
\]
\section{Cauchy problem for the heat equation at polar coordinates.}

\subsection{Auxiliary propositions}

We define polynomials $\mbox{W}_j \left( \mbox{z} \right)$ with the help of
the generating function $e^{-t^2}I_0 (2tz)$
\begin{equation}
\label{eq7}
e^{-t^2}I_0 (2tz)=\sum\limits_{j=0}^\infty {\frac{t^{2j}}{(2j)!}} W_j (z).
\end{equation}
where $I_0 (x)$ - the zero-order Bessel function of the first kind.

It follows from (\ref{eq7}) that the polynomials $W_j (z)$ have the form
\begin{equation}
\label{eq8}
W_j (z)=\frac{d^{2j}}{dt^{2j}}\left[ {e^{-t^2}I_0 (2tz)} \right]_{t=0}
\end{equation}
We get another view of polynomials $W_j (z)$. From the definition of the
Bessel operator, [1],
\[
B=\frac{d^2}{dz^2}+\frac{1}{z}\frac{d}{dz}
\]
it follows that
\[
B^j[I_0 (2tz)]=(2t)^{2j}I_0 (2tz),
\]
Therefore
\[
\exp \left( {-\frac{B}{4}} \right)[I_0 (2tz)]=e^{-t^2}I_0 (2tz),
\]
In this formula we equate the coefficients of the $\xi ^{2j}$ degree at the
left and right sides. We get
\[
\exp \left( {-\frac{B}{4}} \right)\left[ {\frac{2^{2j}z^{2j}}{2^{2j}j!^2}}
\right]=\frac{W_j (z)}{(2j)!},
\]
hence the polynomials$\mbox{W}_j \left( \mbox{z} \right)$ have the form
\[
W_j (z)=\frac{(2j)!}{j!^2}\exp \left( {-\frac{B}{4}} \right)[z^{2j}].
\]
With the help of the polynomials $\mbox{W}_j \left( \mbox{z} \right)$ we will
obtain the new formulas for solutions of direct and inverse Cauchy problems
at the polar coordinates.

\subsection{New formulas for solution of the Cauchy problem at polar
coordinates}

We deduce the new formulas for the solution of the Cauchy problem for the
heat equation at polar coordinates $\left( {r,\phi } \right)$ if the thermal
regime depends only on the variable $r$. We use the explicit formula for the
Cauchy problem solution at polar coordinates [6]
\[
u\left( {\tau ,r} \right)=\int_0^\infty \lambda e^{-\lambda ^2t}J_0 (\lambda
r)\left( {\int_0^\infty \xi J_0 (\lambda \xi )\,f\left( \xi \right)d\xi }
\right)d\lambda ,
\]
where $J_0 (\lambda \xi )$ - the zero-order modified Bessel functions, [1].

We write the last equation in the form
\begin{equation}
\label{eq9}
u\left( {\tau ,r} \right)=\int_0^\infty \lambda e^{-\lambda ^2(t+\beta
)}e^{\lambda ^2\beta }J_0 (\lambda r)\left( {\int_0^\infty \xi J_0 (\lambda
\xi )\,f\left( \xi \right)d\xi } \right)d\lambda ,
\end{equation}
where $\beta >0$.

In (\ref{eq7}) we make the substitutions: $t=i\lambda \sqrt \beta $ ,
$z=\frac{x}{2\sqrt \beta }$. Then we obtain
\begin{equation}
\label{eq10}
e^{\lambda ^2\beta }J_0 (\lambda x)=\sum\limits_{j=0}^\infty (
-1)^j\frac{\lambda ^{2j}\beta ^j}{(2j)!}W_j \left( {\frac{x}{2\sqrt \beta }}
\right).
\end{equation}
The formula (\ref{eq9}) according to the relation (\ref{eq10}) takes the form
\begin{equation}
\label{eq11}
u\left( {\tau ,r} \right)=\int_0^\infty {e^{-\lambda ^2(\tau +\beta )}} J_0
(\lambda x)\sum\limits_{j=0}^\infty ( -1)^j\frac{\lambda ^{2j}}{(2j)!}\beta
^j\int_0^\infty {W_j } \left( {\frac{\xi }{2\sqrt \beta }} \right)f\left(
\xi \right)d\xi d\lambda
\end{equation}
At the last formula change the integration order and compute the inner
integral on the variable $\lambda $. The Weber integral [6] is used in
calculation.
\[
\int_0^\infty \lambda e^{-\lambda ^2(\tau +\beta )}J_0 (\lambda r)J_0
(\lambda \xi )d\lambda =\frac{e^{-\frac{r^2+\xi ^2}{4(\tau +\beta
)}}}{2(\tau +\beta )}I_0 \left( {\frac{r\xi }{2(\tau +\beta )}} \right).
\]
By equating the coefficients before the $\xi ^{2j}$ degree we get the value
of an integral
\begin{equation}
\label{eq12}
\frac{1}{j!^22^{2j}}\int_{-\infty }^\infty ( -1)^je^{-\lambda ^2(\tau +\beta
)}J_0 (\lambda r)\lambda ^{2j}d\lambda =\frac{1}{(2j)!}\frac{d^{2j}}{d\xi
^{2j}}\left[ {\frac{e^{-\frac{r^2+\xi ^2}{4(\tau +\beta )}}}{2(\tau +\beta
)}I_0 \left( {\frac{r\xi }{2(\tau +\beta )}} \right)} \right]_{\xi =0} .
\end{equation}
Due to ~formula (\ref{eq8}) the last formula can be written
\[
\frac{e^{-\frac{r^2}{4(\tau +\beta )}}}{2(\tau +\beta )}\frac{1}{2^{2j}(\tau
+\beta )^j(2j)!}W_j \left( {\frac{r}{2\sqrt {\tau +\beta } }}
\right)=\frac{1}{(2j)!}\frac{d^{2j}}{d\xi ^{2j}}\left[
{\frac{e^{-\frac{r^2+\xi ^2}{4(\tau +\beta )}}}{2(\tau +\beta )}I_0 \left(
{\frac{r\xi }{2(\tau +\beta )}} \right)} \right]_{\xi =0}
\]
Then (\ref{eq12}) becomes
\begin{equation}
\label{eq13}
\int_{-\infty }^\infty ( -1)^je^{-\lambda ^2(\tau +\beta )}J_0 (\lambda
r)\lambda ^{2j}d\lambda =\frac{e^{-\frac{r^2}{4(\tau +\beta )}}}{2(\tau
+\beta )}\frac{j!^2}{(\tau +\beta )^j(2j)!}W_j \left( {\frac{r}{2\sqrt {\tau
+\beta } }} \right).
\end{equation}
Finally, formula (\ref{eq11}) for the thermal field $u(\tau ,r)$ with (\ref{eq13}) takes
the form
\begin{equation}
\label{eq14}
u(\tau ,r)=\frac{e^{-\frac{r^2}{4(\tau +\beta )}}}{2(\tau +\beta
)}\sum\limits_{j=0}^\infty {\frac{1}{(\tau +\beta )^j}} W_j \left(
{\frac{r}{2\sqrt {\tau +\beta } }} \right)\frac{j!^2\beta ^j}{(2j)!^2}f_j ,
\end{equation}
where
\[
f_j =\int_{-\infty }^\infty {W_j } \left( {\frac{\xi }{2\sqrt \beta }}
\right)f\left( \xi \right)d\xi .
\]
\textbf{Remark}. From the last equality for $u(\tau ,r)$ at $\tau =0$ the
expansion theorem on the eigenfunctions $\left\{ {e^{-\frac{r^2}{4\beta
}}W_j \left( {\frac{r}{2\sqrt \beta }} \right)} \right\}$ can be obtained
\[
f(r)=\frac{e^{-\frac{r^2}{4\beta }}}{2\beta }\sum\limits_{j=0}^\infty {W_j }
\left( {\frac{r}{2\sqrt \beta }} \right)\frac{j!^2}{(2j)!^2}f_j ,
\]
where
\[
f_j =\int_{-\infty }^\infty {W_j } \left( {\frac{\xi }{2\sqrt \beta }}
\right)f\left( \xi \right)d\xi .
\]
Similarly to section-1 we can get two formulas for the thermal field. If we
replace in the formula (\ref{eq14}): $\beta \leftrightarrow \tau +\beta $, then the
new formula takes the form
\begin{equation}
\label{eq15}
u(\tau ,r)=\frac{e^{-\frac{r^2}{4\beta }}}{2\beta }\sum\limits_{j=0}^\infty
{\frac{1}{\beta ^j}} W_j \left( {\frac{r}{2\sqrt \beta }}
\right)\frac{j!^2\left( {\tau +\beta } \right)^j}{(2j)!^2}f_j ,
\end{equation}
where
\[
f_j =\int_{-\infty }^\infty {W_j } \left( {\frac{\xi }{2\sqrt {\tau +\beta }
}} \right)f\left( \xi \right)d\xi .
\]
The third new formula is proved similarly to the section-1. We apply the
formula from [10]
\[
J_0 (\lambda x)J_0 (\lambda y)=\frac{1}{\pi }\int_0^\pi {J_0 (\lambda \sqrt
{x^2+y^2-2xy\cos \phi } d\phi }
\]
to the formula (\ref{eq9}).

Then
\begin{equation}
\label{eq16}
e^{\lambda ^2\tau }J_0 (\lambda x)J_0 (\lambda y)=\frac{1}{\pi
}\sum\limits_{j=0}^\infty ( -1)^j\frac{\lambda ^{2j}\tau
^j}{(2j)!}\int_0^\pi {W_j } \left( {\frac{\sqrt {x^2+y^2-2xy\cos \phi }
}{2\sqrt \tau }} \right)d\phi .
\end{equation}
Write equation (\ref{eq9}) as
\[
u\left( {\tau ,r} \right)=\int_0^\infty \lambda e^{-\lambda ^2(t+\beta
)}\left( {\int_0^\infty {e^{\lambda ^2\beta }} J_0 (\lambda r)J_0 (\lambda
\xi )\,\xi f\left( \xi \right)d\xi } \right)d\lambda ,
\]
where $\beta >0$.

Next, on the basis of (\ref{eq16}) we have
\[
u\left( {\tau ,r} \right)=\frac{1}{\pi }\int_0^\infty {e^{-\lambda ^2(\tau
+\beta )}} \sum\limits_{j=0}^\infty ( -1)^j\frac{\lambda ^{2j}}{(2j)!}\beta
^j\int_0^\infty {\int_0^\pi {W_j } } \left( {\frac{\sqrt {x^2+y^2-2xy\cos
\phi } }{2\sqrt \beta }} \right)d\phi \xi f\left( \xi \right)d\xi d\lambda
\]
We change the order of integration and use the definition of the
Gamma-function [10] to compute the inner integral on the variable $\lambda .$
\[
\int_0^\infty ( -1)^je^{-\lambda ^2(\tau +\beta )}\lambda ^{2j}d\lambda
=\frac{\left( {-1} \right)^j\Gamma \left( {j+\frac{1}{2}} \right)}{2(\tau
+\beta )^{j+\frac{1}{2}}}.
\]
Finally, we obtain
\begin{equation}
\label{eq17}
u(\tau ,r)=\frac{1}{2}\sum\limits_{j=0}^\infty {\left( {-1}
\right)^j\frac{\beta ^j}{(\tau +\beta )^{j+\frac{1}{2}}}} \frac{\Gamma
\left( {j+\frac{1}{2}} \right)}{(2j)!}f_j ,
\end{equation}
Where
\[
f_j =\int_0^\infty {\int_0^\pi {W_j } } \left( {\frac{\sqrt {r^2+\xi
^2-2r\xi \cos \phi } }{2\sqrt \beta }} \right)d\phi \xi f\left( \xi
\right)d\xi .
\]
\section{The inverse Cauchy problem for the heat equation at polar
coordinates.}

The inverse Cauchy problem [3],[4],[5] in polar coordinates leads to the solution of the
first- kind integral Fredholm equation
\[
\int_0^\infty {\frac{e^{-\frac{r^2+\xi ^2}{4\tau }}}{2\tau }} I_0 \left(
{\frac{r\xi }{2\tau }} \right)f(\xi )d\xi =u(\tau ,\xi ).
\]
To find three previously unknown expressions for the solution of inverse
Cauchy problem in the series on the polynomials $W_j (z)$(\ref{eq8}) we will do
the same in section-1.

The first expression is
\begin{equation}
\label{eq18}
f(r)=\frac{e^{-\frac{r^2}{4\beta }}}{2\beta }\sum\limits_{j=0}^\infty
{\frac{1}{\beta ^j}} W_j \left( {\frac{r}{2\sqrt \beta }}
\right)\frac{j!^2(\beta +\tau )^j}{(2j)!^2}u_j ,
\end{equation}
where
\[
u_j =\int_0^\infty {W_j } \left( {\frac{\xi }{2\sqrt {\tau +\beta } }}
\right)u\left( {\tau ,\xi } \right)d\xi .
\]
The second expression is
\begin{equation}
\label{eq19}
f(r)=\frac{e^{-\frac{r^2}{4(\beta +\tau )}}}{2(\beta +\tau
)}\sum\limits_{j=0}^\infty {\frac{1}{(\beta +\tau )^j}} W_j \left(
{\frac{r}{2\sqrt {(\beta +\tau )} }} \right)\frac{j!^2\tau ^j}{(2j)!^2}u_j
,
\end{equation}
where
\[
u_j =\int_0^\infty {W_j } \left( {\frac{\xi }{2\sqrt \beta }} \right)u\left(
{\tau ,\xi } \right)d\xi .
\]
The third expression is
\begin{equation}
\label{eq20}
f(r)=\frac{1}{2}\sum\limits_{j=0}^\infty {\left( {-1} \right)^j\frac{\left(
{\tau +\beta } \right)^j}{\tau ^{j+\frac{1}{2}}}} \frac{\Gamma \left(
{j+\frac{1}{2}} \right)}{(2j)!}u_j ,
\end{equation}
where
\[
u_j =\int_0^\infty {\int_0^\pi {W_j } } \left( {\frac{\sqrt {r^2+\xi
^2-2r\xi \cos \phi } }{2\sqrt {\tau +\beta } }} \right)d\phi \xi u\left(
{\tau ,\xi } \right)d\xi .
\]

\end{document}